\documentstyle [12pt,epsf, graphicx, amssymb,hyperref,todonotes, amsmath]{article}

\begin{document}
\def\l{\lambda}
\def\m{\mu}
\def\a{\alpha}
\def\b{\beta}
\def\g{\gamma}
\def\d{\delta}
\def\e{\epsilon}
\def\o{\omega}
\def\O{\Omega}
\def\v{\varphi}
\def\t{\theta}
\def\r{\rho}
\def\bs{$\blacksquare$}
\def\bp{\begin{proposition}}
\def\ep{\end{proposition}}
\def\bt{\begin{th}}
\def\et{\end{th}}
\def\be{\begin{equation}}
\def\ee{\end{equation}}
\def\bl{\begin{lemma}}
\def\el{\end{lemma}}
\def\bc{\begin{corollary}}
\def\ec{\end{corollary}}
\def\pr{\noindent{\bf Proof: }}
\def\note{\noindent{\bf Note. }}
\def\bd{\begin{definition}}
\def\ed{\end{definition}}
\def\C{{\mathbb C}}
\def\P{{\mathbb P}}
\def\Z{{\mathbb Z}}
\def\d{{\rm d}}
\def\deg{{\rm deg\,}}
\def\deg{{\rm deg\,}}
\def\arg{{\rm arg\,}}
\def\min{{\rm min\,}}
\def\max{{\rm max\,}}

\newcommand{\norm}[1]{\left\Vert#1\right\Vert}
\newcommand{\abs}[1]{\left\vert#1\right\vert}

\newcommand{\set}[1]{\left\{#1\right\}}
\newcommand{\setb}[2]{ \left\{#1 \ \Big| \ #2 \right\} }

\newcommand{\IP}[1]{\left<#1\right>}
\newcommand{\Bracket}[1]{\left[#1\right]}
\newcommand{\Soger}[1]{\left(#1\right)}

\newcommand{\Integer}{\mathbb{Z}}
\newcommand{\Rational}{\mathbb{Q}}
\newcommand{\Real}{\mathbb{R}}
\newcommand{\Complex}{\mathbb{C}}

\newcommand{\eps}{\varepsilon}
\newcommand{\To}{\longrightarrow}
\newcommand{\varchi}{\raisebox{2pt}{$\chi$}}

\newcommand{\E}{\mathbf{E}}
\newcommand{\Var}{\mathrm{var}}

\def\squareforqed{\hbox{\rlap{$\sqcap$}$\sqcup$}}
\def\qed{\ifmmode\squareforqed\else{\unskip\nobreak\hfil
\penalty50\hskip1em\null\nobreak\hfil\squareforqed
\parfillskip=0pt\finalhyphendemerits=0\endgraf}\fi}

\renewcommand{\th}{^{\mathrm{th}}}
\newcommand{\Dif}{\mathrm{D_{if}}}
\newcommand{\Difp}{\mathrm{D^p_{if}}}
\newcommand{\GHF}{\mathrm{G_{HF}}}
\newcommand{\GHFP}{\mathrm{G^p_{HF}}}
\newcommand{\f}{\mathrm{f}}
\newcommand{\fgh}{\mathrm{f_{gh}}}
\newcommand{\T}{\mathrm{T}}
\newcommand{\K}{^\mathrm{K}}
\newcommand{\PghK}{\mathrm{P^K_{f_{gh}}}}
\newcommand{\Dig}{\mathrm{D_{ig}}}
\newcommand{\for}{\mathrm{for}}
\newcommand{\End}{\mathrm{end}}

\newtheorem{th}{Theorem}[section]
\newtheorem{lemma}{Lemma}[section]
\newtheorem{definition}{Definition}[section]
\newtheorem{corollary}{Corollary}[section]
\newtheorem{proposition}{Proposition}[section]

\begin{titlepage}

\begin{center}

\topskip 5mm

{\LARGE{\bf {Higher derivatives of functions

\vskip 4mm

with given critical points and values}}}

\vskip 8mm

{\large {\bf G. Goldman}}

\vspace{6 mm}

{Department of Applied Mathematics, Tel Aviv University,
Tel Aviv 69978, Israel. e-mail: ggoldman@tauex.tau.ac.il}

\vspace{6 mm}

{\large {\bf Y. Yomdin}}

\vspace{6 mm}

{Department of Mathematics, The Weizmann Institute of Science,
Rehovot 76100, Israel. e-mail: yosef.yomdin@weizmann.ac.il}

\end{center}

\vspace{6 mm}
\begin{center}

{ \bf Abstract}
\end{center}

{\small Let $f: B^n \rightarrow {\mathbb R}$ be a $d+1$ times continuously differentiable function on the unit ball $B^n$, with $\max_{z\in B^n} \Vert f(z) \Vert=1$. A well-known fact is that if $f$ vanishes on a set $Z\subset B^n$ with a non-empty interior, then for each $k=1,\ldots,d+1$ the norm of the $k$-th derivative $||f^{(k)}||$ is at least $M=M(n,k)>0$. A natural question to ask is ``what happens for other sets $Z$?''. This question was partially answered in \cite{Yom5}-\cite{Yom7}). In the present paper we ask for a similar (and closely related) question: what happens with the high-order derivatives of $f$, if its gradient vanishes on a given set $\Sigma$? And what conclusions for the high-order derivatives of $f$ can be obtained from the analysis of the metric geometry of the ``critical values set'' $f(\Sigma)$? In the present paper we provide some initial answers to these questions}.

\end{titlepage}

\newpage


\section{Introduction}\label{Sec:Intro}
\setcounter{equation}{0}


\smallskip


\smallskip

In this paper we continue our study (started in \cite{Yom1},\cite{Yom5}-\cite{Yom7}) of certain very special settings of the classical Whitney's smooth extension problem (see \cite{Bru.Shv,Fef,Fef.Kla,Whi1,Whi2,Whi3}). Let's recall shortly the general setting of this problem, and its current status.

\medskip

Let $E\subset B^n \subset {\mathbb R}^n$ be a closed subset of the unit ball $B^n$, and let $\bar f$ be a real function, defined on $E$. Can $\bar f$ be extended to a $C^d$-smooth $f$ on ${\mathbb R}^n$, and, if the extension is possible, what is the minimal $C^d$-norm of $f$?

\medskip

Recent exciting developments in the general Whitney problem (see \cite{Bru.Shv,Fef,Fef.Kla} and references therein), provide essentially a complete answer to the general Whitney extension problem in any dimension. In particular, the results of \cite{Bru.Shv,Fef,Fef.Kla}), provide an important information on this problem, which was earlier available only in dimension one. The "finiteness principle" achieved in this recent work, claims that, as in classical Whitney's results in dimension one (\cite{Whi2}), it is enough to check only finite subsets of $Z$ with cardinality bounded in terms of $n$ and $d$ only. There is also an algorithmic way to estimate the minimal extension norm for any finite $Z$.

\medskip

However, a possibility of an explicit answer, as in dimension one, through a kind of multi-dimensional divided finite differences, remains an open problem.

\smallskip
\medskip

Let us now describe the setting of \cite{Yom5}-\cite{Yom7}, which we use below. Let $Z\subset B^n \subset {\mathbb R}^n$ be a closed subset of the unit ball $B^n$. In \cite{Yom5}-\cite{Yom7} we look for $C^{d+1}$-smooth functions $f:B^n\to {\mathbb R}$, vanishing on $Z$. Such $C^{d+1}$-smooth (and even $C^\infty$) functions $f$ always exist, since any closed set $Z$ is a set of zeroes of a $C^\infty$-smooth function.

\medskip

We normalize the extensions $f$ requiring $\max_{B^n}|f|=1$, and ask for the minimal possible norm of the last derivative $||f^{(d+1)}||$, which we call {\it the $d$-rigidity ${\cal RG}_d(Z)$ of $Z$}. In other words, for each normalized $C^{d+1}$-smooth function $f:B^n\to {\mathbb R},$ vanishing on $Z$, we have
$$
||f^{(d+1)}||\ge {\cal RG}_d(Z),
$$
and ${\cal RG}_d(Z)$ is the maximal number with this property.

\medskip

Our previous papers \cite{Yom5,Yom6,Yom7}, related to smooth rigidity, provide certain bounds on ${\cal RG}_d(Z)$, in terms of the fractal geometry of $Z$. We use some of these results below. In order to compare the results of \cite{Yom5,Yom6,Yom7} with the general results, available today in Whitney's extension theory, let's make the following remark:

\medskip

Of course, the results of \cite{Fef,Fef.Kla} provide, in principle, an algorithmic way to estimate also our quantities ${\cal RG}_d(Z)$, for any closed $Z\subset B^n$ (via considering finite subsets of $Z$ of bounded cardinality). However, our goal in the present paper, as well as in our previous papers \cite{Yom1,Yom5,Yom6,Yom7}, related to smooth rigidity, is somewhat different: {\it we look for an explicit answer in terms of simple, and directly computable geometric (or topological) characteristics of $Z$}.

\medskip

Now we finally come to the setting of the problem in the present paper, and to the new results below. In the present paper we consider, as above, $C^{d+1}$-smooth functions $f:B^n\to {\mathbb R}$. But, in contrast to \cite{Yom5,Yom6,Yom7}, we do not consider zero sets of $f$. Instead, we assume some geometric conditions on the critical points and values of $f$, and, as above, derive in conclusion some lower bounds on $||f^{(d+1)}||$.

\medskip

More accurately, put, as usual,
$$
\nabla f(x)=(\frac {\partial f}{\partial x_1}(x),\ldots,\frac {\partial f}{\partial x_n}(x))
$$
at each point $x\in B^n$.

\medskip

The point $x\in B^n$ is called a critical (or a singular) point for $f$, if the vector equation $\nabla f(x)=0$ is satisfied. 
The real $\nu$ is called a critical (or a singular) value of $f$ if $\nu=f(x)$ for a certain critical point $x$ of $f$.

\medskip

\noindent {\bf Remark} {\it Prescribing the set of the critical points of $f$ is not that immediate as for the zeroes set of $f$ (which can be any closed subset of $B^n$). Indeed, consider, for instance $Z=S=\{x_1^2+x_2^2=\frac{1}{4}\}$ being the circle of radius $\frac{1}{2}$. Then for any smooth $f$ on $B^2$ with $\nabla f$ vanishing on $S$ the restriction of $f$ to $S$ is a constant. Consequently, there is a critical point of $f$ inside $S$.

\medskip

For $Z=T^2 \subset B^3$ being the standard torus, any smooth $f$ on $B^2$ with $\nabla f$ vanishing on $Z$ must have some critical points in the interior of the torus, by similar topological reasons.}

\medskip


\medskip
\medskip

Let us now state the main results of this paper. From now on we always assume that an integer $d\ge 1$ is fixed. We consider two settings of the problem:

\smallskip

First, in Section \ref{Sec:crit.pts}, we look only on the geometry of the critical points of $f$, not taking into account the critical values, and not requiring any non-degeneracy. Here we can apply, with almost no changes, the results of \cite{Yom5}-\cite{Yom7}, and obtain a series of new results, where the set $Z$ of zeroes of $f$ is replaced by the set $\Sigma$ of the critical points of $f$.

\medskip


\smallskip

Next, in Section \ref{Sec:Crit.values}, we look only on the geometry of the critical values of $f$, completely ignoring the geometry of the critical points. Here, the results are obtained by a kind of "backward reading" of the old results of \cite{Yom01,Yom.Com}. Still, these new results seems to provide an important new information on the rigidity properties of smooth functions.

\smallskip


\section{Geometry of critical points}\label{Sec:crit.pts}
\setcounter{equation}{0}

Let $\Sigma \subset B^n \subset {\mathbb R}^n$ be a closed subset of the unit ball $B^n$. We look for $C^{d+1}$-smooth functions $f:B^n\to {\mathbb R}$, with the gradient $f$ vanishing on $\Sigma$.

\medskip

We normalize these functions $f$ requiring $\max_{B^n}|f|=1$, on one side, and $f(x_0)=0$ for certain $x_0\in B^n$, on the other side, and call the set of such functions $U(d, \Sigma)$. Then we ask for the minimal possible norm of the last derivative $||f^{(d+1)}||$, which we call {\it the $d^1$-rigidity ${\cal RG}^1_d(\Sigma)$ of $\Sigma$}.

\medskip

To avoid possible misreading, let's provide a formal definition:

\bd\label{def:rig1}
The $d^1$-rigidity ${\cal RG}^1_d(\Sigma)$ of $\Sigma$ is defined as
$$
{\cal RG}^1_{d-1}(\Sigma) = \max R, with \ \ ||f^{(d+1)}|| \ge R, \ \  {for \ each} \ \ f \in U(d, \Sigma).
$$
\ed

In other words, for each normalized $C^{d+1}$-smooth function $f:B^n\to {\mathbb R},$ with $\nabla f$ vanishing on $\Sigma$, we have
$$
||f^{(d+1)}||\ge {\cal RG}^1_d(\Sigma).
$$
Now we show that the ``rigidities'' ${\cal RG}^1_d(\Sigma)$ and ${\cal RG}_d(\Sigma)$ are ``subordinated'': the second majorates the first:

\bt\label{Thm.General.Deriv}
$$
{\cal RG}^1_{d}(\Sigma)\ge \frac{1}{\sqrt n} {\cal RG}_d(\Sigma).
$$
\et
\pr
Pick a function $f \in U(d, \Sigma)$. Now we notice that if $\nabla f$
vanishes on $\Sigma$, then each of the partial derivatives $\frac {\partial f}{\partial x_i}$ vanishes on $\Sigma$. On the other hand, for the normalized $f$, at least one of these partial derivatives, say, $f'_1=\frac {\partial f}{\partial x_1}$, attains sufficiently big value, say, $\frac{1}{\sqrt n}$, inside the ball $B^n$. Indeed, when we integrate from $x_0$ (where $f$ vanishes), to $x_1$, where $|f(x_1)|=1$, the integral, which is the sum of the integrals of the partial derivatives, with the weights $\le 1$, is equal to one.

\medskip

Therefore, the function $f'_1$ vanishes on $\Sigma$, and its maximum on $B^n$ is at least $\frac{1}{\sqrt n}$. According to the definition of ${\cal RG}_d(\Sigma)$ given above, we conclude that
$$
||(f'_1)^{(d)}|| \ge \frac{1}{\sqrt n} {\cal RG}_d(\Sigma),
$$
which implies
$$
{\cal RG}^1_{d}(\Sigma)\ge \frac{1}{\sqrt n} {\cal RG}_d(\Sigma).
$$
This completes the proof of Theorem \ref{Thm.General.Deriv}. $\square$

\medskip

Our main conclusion is that the results of \cite{Yom5}-\cite{Yom7} are directly applicable to the partial derivatives $\frac {\partial f}{\partial x}$. These results are given in terms of the metric ``density'' of $\Sigma$, of its ``Remez Constant'', its metric entropy (the asymptotic behavior of the covering numbers), and in terms of the topology of $\Sigma$. We address the interested reader to \cite{Yom5}-\cite{Yom7}.


\smallskip


\section{Critical values}\label{Sec:Crit.values}
\setcounter{equation}{0}

This section presents the main new results of the present paper. Because of the nature of the results and of the proofs, we extend the original setting above from functions $f$ on $B^n$ to mappings $f:B^n_r\to B^m, \ m\le n, \ r>0.$

\medskip

For a mapping $f$ as above, and for each $x\in B_r^n$, the differential $df(x)$ is a linear mapping $df(x): {\mathbb R}^n\to {\mathbb R}^m.$ The image $df(x)(B^n)$ of the unit ball $B^n$ in ${\mathbb R}^n$ is an ellipsoid in ${\mathbb R}^m$. We denote the semi-axes of this ellipsoid by

$$
\lambda_0:=1, \ \lambda_1(x)\le \lambda_2(x)\le \ldots \le \lambda_m(x),
$$
and put
$$
\Lambda(f,x)= \ (\lambda_1(x),\lambda_2(x), \ldots, \lambda_m(x)).
$$

Accordingly, we define the set of near-critical points and values of $f$ as follows:

\bd\label{def:crit}

For a given $\Lambda=(\lambda_1,\lambda_2, \ldots, \lambda_m)$ the set $\Sigma(f,\Lambda)$ of $\Lambda$-near-critical points of $f$, and the set $\Delta(f,\Lambda)$ of $\Lambda$-near-critical values of $f$ are defined as follows:

$$
\Sigma(f,\Lambda) = \{x\in B_r^n, \ \lambda_i(x) \le \lambda_i, \ i=1,\ldots,m \}
$$

\medskip

$$
\Delta(f,\Lambda)=f(\Sigma(f,\Lambda)).
$$
\ed

Before stating the main results of \cite{Yom01,Yom.Com}, lets recall the definition of two of the most important constants in our further calculations:

\bd\label{def:taylor.const}
For a natural $d$ the Taylor constant $R_d(f)$ is defined as
$$
R_d(f)=\frac{\max_{x\in B_r^n} ||f^{(d)}||}{d!}r^d.
$$
\ed

\bd\label{def:covering.number} For a compact set $S$ in a metric space $Y$, and for $\e>0$ the covering number $M(\e,S)$ is defined as the minimal number of the $\e$-balls in $Y$, covering $S$.
\ed

Now we recall the main result of \cite{Yom01,Yom.Com} which is crucial for our results below. For a convenience reasons we always put $\lambda_0=1$.

\bt\label{thm.old.main} (Theorem 9.2, \cite{Yom.Com})

\smallskip

Let $f:B^n_r\to {\mathbb R}^m$ be a $C^d$-smooth mapping. Then for a given
$$
\Lambda=(\lambda_1,\lambda_2, \ldots, \lambda_m)
$$
and for $\e>0$ we have

$$
M(\e, \Delta(f,\Lambda)) \le c\sum_{i=0}^m \lambda_0\lambda_1\lambda_2 \ldots \lambda_i(\frac{r}{\e})^i, \ \ \ \ \ \ \ \e \ge R_d(f).
$$

\medskip

$$
M(\e, \Delta(f,\Lambda)) \le c\sum_{i=0}^m \lambda_0\lambda_1\lambda_2 \ldots \lambda_i(\frac{r}{\e})^i
(\frac{R_d(f)}{\e})^{\frac{n-i}{d}}, \ \ \ \ \ \ \ \e \le R_d(f),
$$
where $c$ is a constant depending only on $n,d$. 
\et
Our goal now is ``to read Theorem \ref{thm.old.main} in the opposite direction'', and to bound $R_d(f)$ from below in terms of
$M(\e, \Delta(f,\Lambda))$.

To achieve this goal, first of all we notice that we are interested only in the case $\e \le R_d(f)$. Indeed, for $\e \ge R_d(f)$ the first inequality of Theorem \ref{thm.old.main} produces the bound on $M(\e, \Delta(f,\Lambda))$, which is true for $f$ being a polynomial of degree $d-1$ (See \cite{Yom.Com}). Consequently, in this case we cannot produce any lower bound for $R_d(f)$ from below in terms of
$M(\e, \Delta(f,\Lambda))$.

\medskip

Thus we define the set $E=E(f,\Lambda)) \subset {\mathbb R}_+$ consisting of all $\e$ for which
\be\label{eq:first}
M(\e, \Delta(f,\Lambda)) > c\sum_{i=0}^m \lambda_0\lambda_1\lambda_2 \ldots \lambda_i(\frac{r}{\e})^i.
\ee
In particular, for each $\e \in E$ we have $\e \le R_d(f)$, because otherwise the first inequality of Theorem \ref{thm.old.main} would be violated.

\medskip

Next we notice that the right-hand side of the second inequality of Theorem \ref{thm.old.main} is a polynomial in $\eta: =  (\frac{(R_d(f))}{\e})^{\frac{1}{d}}$ with positive coefficients, and hence it is monotone in $\eta$. Therefore, for any given left-hand side $\nu(\e):=M(\e, \Delta(f,\Lambda))$  of this inequality we get at most one solution $\eta=\eta(\e)$ of the equation
\be\label{eq:main}
\nu(\e) = c\sum_{i=0}^m \lambda_0\lambda_1\lambda_2 \ldots \lambda_i(\frac{r}{\e})^i\eta^{\frac{n-i}{d}}.
\ee
In fact, for $\eta=1$ the right hand side of (\ref{eq:main}) is equal to the right hand side of (\ref{eq:first}), and therefore for each $\e \in E$ we obtain the unique solution $\eta=\eta(\e) > 1$ of (\ref{eq:main}). Of course, the condition $\eta > 1$ is the same as $\e \le R_d(f)$.

\medskip

Now we come to our main result: we assume that the covering number of the critical values of $f$, i.e. $\nu(\e):=M(\e, \Delta(f,\Lambda))$, is known for each $\e>0$, and we would like to extract from this information the lower bound for $\eta: = R_d(f)$.

\bt\label{thm.main}
Let $f$ and $\Lambda$ as above be given. Assume that the set $E=E(f,\Lambda))$ is not empty, and define $\gamma=\gamma(f,\Lambda))$ as
$$
\gamma= \max_{\e \in E} \ \ \  \eta(\e)\e.
$$
Then we have

$$
R_k(f)\ge \gamma.
$$
\et
\pr
The proof of Theorem \ref{thm.main} follows directly from the presentation before the theorem. $\square$

\subsection{Some examples}

The structure of the critical and near-critical values of smooth functions may be pretty complicated (see, e.g. \cite{Yom.Com} and references therein). Accordingly, we consider the problem of understanding the fractal geometry of $\Delta(f)$, via Theorem \ref{thm.main}, from the point of view of smooth rigidity, as an important and non-trivial one. In this paper we provide some initial examples, which may be instructive.

\medskip

First we restrict consideration to the case $m=1$, i.e. to smooth functions $f:B^n_1\to {\mathbb R}$. Second, we consider critical values of $f$ (and not near-critical ones). That is, we put $\lambda_1=0$. Theorem \ref{thm.old.main} takes now the form:

\begin{align}\label{Eq:spec.case}
    M(\e,\Delta(f)) \le c(n,d), &  & \e \ge R_d(f),\\    
    M(\e,\Delta(f)) \le c(n,d) (\frac{(R_d(f))}{\e})^{\frac{n}{d}}, & & \e \le R_d(f).
\end{align}


This is the time to discuss shortly the constant $c(n,d)$. Some ``explicit'' bounds for $c(n,d)$ are given in \cite{Yom.Com}. However, these bounds involve certain rather cumbersome expressions. We consider producing more instructive bounds in the results of \cite{Yom.Com} as an important open problem (which recently attracted some attention in the field). This is the most important reason for us to consider below the case $n=1$, where $c(1,d)=d+1$ is explicit.

\medskip

Now we proceed as follows: we consider a certain subset $\Delta \subset {\mathbb R}$, and find (or estimate) the metric entropy $M(\e,\Delta)$ for various $\e$. On this base we find the set $E$, and apply Theorem \ref{thm.main}. The conclusion is as follows:

\medskip

\noindent  {\it For each smooth $f$ with $\Delta(f)) = \Delta$ we have $R_d(f)\ge \gamma$.}

\medskip

Assume now that the cardinality of $\Delta$ is strictly greater than $c(n,d)$, and denote by $\e_0=\e(\Delta,n,d)$ the biggest $\e$ for which $M(\e,\Delta)\ge c(n,d)+1$. Notice, that if the cardinality of $\Delta$ is exactly $c(n,d)+1$, then $\e_0$ is the minimal distance between the points of $\Delta$.

\medskip

Next for $\e_0$ as above put

\be\label{eq:prep1}
\gamma_0=(1+\frac{1}{c(n,d)})^{\frac{k}{n}}\e_0.
\ee
Now we can present our first specific corollary of general Theorem \ref{thm.main}:

\bc\label{Cor:main1}
For $\Delta$ and $\e_0$ as above, and for each smooth $f$ with

\noindent $\Delta(f)) = \Delta$ we have $R_d(f)\ge \gamma_0>0$.
\ec
\pr
By definition, $\e_0\in E$, and thus, according to Theorem \ref{thm.main}, we can solve the second equation in (\ref{Eq:spec.case}), and denote its solution by $\gamma_0$:

\be\label{Eq:spec.case11}
c(n,d)+1 = M(\e_0,\Delta) = c(n,d) (\frac{(R_d(f))}{\e_0})^{\frac{n}{d}},
\ee
which gives us

$$
(\frac{(R_d(f))}{\e_0})^{\frac{n}{d}}=\frac{c(n,d)+1}{c(n,d)}=1+\frac{1}{c(n,d)},
$$
which finally provides
$$
\gamma_0=\e_0 [1+\frac{1}{c(n,d)}]^{\frac{d}{n}}.
$$
This completes the proof of Corollary \ref{Cor:main1}. $\square$

\medskip

Our next result deals with even more restricted situation: we consider smooth functions of one variable. Still we believe that the results provided by Theorem \ref{thm.main} in this case are new and instructive. The main advantage of this special case is that the constant $c(1,k)=k+1$ is explicit and accurate, and we state the result in a ``closed form'', not referring to the previous definitions.

\medskip

\bc\label{Cor:main2}
Let $|\Delta|\ge d+1$. Denote by $\e_0=\e(\Delta,1,d)$ the biggest $\e$ for which $M(\e,\Delta)\ge d+2$. Next put
\be\label{eq:prep11}
\gamma_0=(1+\frac{1}{d+1})^d\e_0.
\ee
Then for each smooth $f:B^1 \to {\mathbb R}$ with
\noindent $\Delta(f) = \Delta$,
we have $R_d(f)\ge \gamma_0>0$.
\ec
\pr
This is a special case of Corollary \ref{Cor:main1}. $\square$

\medskip

To illustrate our approach, let's provide even more a special case of the situation above:

\bc\label{Cor:main3}
Let $|\Delta|\ge 6$. Denote by $\e_0$ the biggest $\e$ for which $M(\e,\Delta)\ge 6$. Next put
\be\label{eq:prep11}
\gamma_0=(\frac{7}{6})^5\e_0.
\ee
Then for each $5$-smooth $f:B^1 \to {\mathbb R}$ with
\noindent $\Delta(f) = \Delta$ we have

$$
R_5(f)\ge \gamma_0=(\frac{7}{6})^5\e_0>0.
$$
\ec
\pr
This is a special case of Corollary \ref{Cor:main1}. $\square$

\medskip

Finally we provide another very specific example, which, however, illustrates the power of the lower bounds for the high-order derivatives,
provided by Corollary \ref{Cor:main1}. We are in a situation of Corollary \ref{Cor:main1}, i.e. we consider $d$-smooth functions $f: B^n\to {\mathbb R}$, and sets $\Delta$ of their critical values. Lets consider $\Delta=\Delta_\alpha =\{1, 2^\alpha, 3^\alpha, \ldots, m^\alpha, \ldots\}$. We assume that $\alpha < 0$, so the points of $\Delta$ always converge to $0$. For the differences $m^\alpha - (m+1)^\alpha$ we get
$$
m^\alpha -(m+1)^\alpha= m^\alpha(1-(1+\frac{1}{m})^\alpha=m^\alpha(-\frac{\alpha}{m}-...)\asymp -\alpha m^{(\alpha-1)}.
$$
Next we fix a certain sufficiently small value of $\e>0$, in particular, providing $M(\Delta,\e)>c(n,d)$. By the above calculation we have $M(\Delta,\e) \asymp (\frac{\e}{\alpha})^{\frac{1}{\alpha-1}}.$ Recall that $\alpha, \ \alpha-1 < 0$.

\medskip

Finally, by Theorem \ref{thm.main}, we get

$$
M(\Delta,\e) \asymp (\frac{\e}{\alpha})^{\frac{1}{\alpha-1}} \le c(n,d)(\frac{R_d(f)}{\e})^{\frac{n}{d}},
$$
or
$$
R_d(f)\ge \e c(n,d)^{-\frac {d}{n}} (\frac {\e}{\alpha})^{\frac {d}{n(\alpha-1)}}= c(n,d)^{-\frac {d}{n}} (\frac {1}{\alpha})^{\frac {d}{n(\alpha-1)}} \e^{1+\frac {d}{n(\alpha-1)}}.
$$
What is most important to us in this expression is the dependence of the answer in $\e$ - in our setting all the other parameters are fixed. Thus we get
\be\label{eq:lower.bound}
R_d(f)\ge C  \e^{1+\frac {d}{n(\alpha-1)}}.
\ee
The conclusion is that for $\e\to 0$ we have, essentially, two quite different possibilities:

\medskip

\noindent 1. If the power of $\e$ is positive, then, our lower bound (\ref{eq:lower.bound}) for $R_d(f)$ decreases, as $\e\to 0$, and this does not exclude $\Delta$ to be the set of critical values of a $C^d$-smooth $f$. The examples in \cite{Yom.Com} illustrate this case.

\medskip

\noindent 2. If the power of $\e$ is negative, then, our lower bound (\ref{eq:lower.bound}) for $R_d(f)$ increases, as $\e\to 0$. Of course, this excludes a possibility for $\Delta$ to be the set of critical values of a $C^d$-smooth $f$. The examples in \cite{Yom.Com} illustrate also this case, but we believe, that the explicit lower bounds for the derivatives of $f$, given by (\ref{eq:lower.bound}) are new.




\medskip

\bibliographystyle{amsplain}

\end{document}